\def\be{\begin{equation}}
\def\ee{\end{equation}}
\def\ba{\begin{array}}
\def\ea{\end{array}}
\def\t{\tilde}
\def\Rb{{I\!\! R}}
\def\Cb{\ \hbox{\vrule width 0.6pt height 6pt depth 0pt
                      \hskip -3.5 pt} C}
\documentstyle[12pt]{article}
\begin{document}
\parskip=0pt
\parindent=18pt
\baselineskip=20pt
\setcounter{page}{1}
\centerline{\Large\bf Highest Weight Representations}
\vspace{2ex}
\centerline{\Large\bf of Quantum Current Algebras}
\vspace{6ex}
\centerline{\large{\sf Sergio Albeverio$^\star$} ~~~and~~~ {\sf Shao-Ming Fei}
\footnote{\sf Alexander von Humboldt-Stiftung fellow.\\
\hspace{5mm}On leave from Institute of Physics, Chinese Academy of Sciences,
Beijing}}
\vspace{4ex}
\parindent=40pt
{\sf Institute of Mathematics, Ruhr-University Bochum,
D-44780 Bochum, Germany}\par
\parindent=35pt
{\sf $^\star$SFB 237 (Essen-Bochum-D\"usseldorf);
BiBoS (Bielefeld-Bochum);\par
\parindent=40pt
CERFIM Locarno (Switzerland)}\par
\vspace{6.5ex}
\parindent=18pt
\parskip=5pt
\begin{center}
\begin{minipage}{5in}
\vspace{3ex}
\centerline{\large Abstract}
\vspace{4ex}
We study the highest weight and continuous tensor product
representations of q-deformed Lie algebras through the
mappings of a manifold into a locally compact group.
As an example the highest weight representation of the q-deformed
algebra $sl_q(2,\Cb)$ is calculated in detail.
\end{minipage}
\end{center}
\newpage

\section{Introduction}

The structure and the representation theory of infinite-dimensional Lie
algebras play important roles in mathematics and physics, especially in
quantum field theory, and have been
investigated for several decades by mathematicians and theoretical physicists.
For instance, by generalizing Cartan's classification of simple Lie algebras,
the Kac-Moody algebras can be built (by generators and relations) from a
Cartan
matrix, and possess a very rich representation theory. Among them is
the class of affine Lie algebras, which can be realized as (untwisted or
twisted) loop algebras, i.e., mappings from the circle $S^1$ into semisimple
Lie algebras. For affine Lie algebras a large class of unitary highest weight
representations can be classified, see e.g. \cite{abook}.

Representations of q-deformed affine
Lie algebras, so called quantum affine Lie algebras, have also been studied
recently [2-4].
However the current algebra (also called gauge algebra)
of the mappings from
a manifold $M$, with dim$(M)>1$, into Lie algebras are no longer
Kac-Moody algebras (they do not possess a Cartan matrix). The motivation
of this paper is to
discuss the representation theory of such kind of
q-deformed current algebras\footnote{
we call them quantum current algebras, although the q-deformation and
physical quantization are different in principle \cite{flato}\cite{fg}}.

Let $g$ be a Lie algebra of a locally compact group $G$. The current
algebra with respect to $g$ is an
algebra of mappings $g^X=Map(X,g)$ from some manifold $X$ into $g$.
The representation of this current algebra can be obtained from the
representation of the related current group which are defined to be group of
mappings $G^X=Map(X,G)$ with its group structure provided by
the pointwise multiplication law. For the representation of the current group,
the continuous unitary representations in Hilbert space
$\pi:G^X\to{u}(H)$ are of significance, where ${u}(H)$ denotes the
unitary operators on a real Hilbert space $H$. Such representations provide a
natural and nontrivial non-commutative generalization of the theory of
distributions on $X$\cite{abook}. Moreover, they correspond to the
representations of the
current algebra by self-adjoint operators, via the exponential mapping.
That is, the representations of current algebras can be obtained as the
infinitesimal representations of current groups.

One approach to obtain
continuous unitary representations of current groups is to look at
the continuous
tensor products of representations of $G$. Another approach
is given by the highest weight representation theory\cite{abook}.
For the group $G=SU(n,1)$, and hence for the algebra $g=su(n,1)$,
the connection between these two kinds of representations
has been made by Torr\'esani \cite{torresani}.

A difficulty concerning the extension of these concepts and relations to
the case of q-deformed Lie algebras and groups lies in that
the representations of the q-deformed Lie algebra can not be obtained
as infinitesimal representations of quantum groups, in fact,
the exponential mappings of
the q-deformed Lie algebra do not give rise, in general, to
representations of quantum groups. It is however
interesting to ask the question whether the quantum current algebras have
also continuous tensor product representations. In this letter we
investigate the construction of
continuous tensor product representations and show that in a particular case
the connection can be made between continuous tensor product representations
and highest weight representations of quantum current algebras.

\section{The continuous tensor product representation of q-deformed current
algebras}

Let ${A}$ be a Lie algebra and ${A}_q$ the q-deformed Lie algebra
of ${A}$ in the sense of \cite{jimbo}\cite{qg}.
We denote by ${B}$ the Borel subalgebra of ${A}_q$
with an involutive antilinear antiautomorphism $\omega$ such that
${B}+\omega.{B}=A$.

{\sf [Definition 1]}. Let $\Lambda$ be a complex character of ${B}$,
and $\pi:A_q\to End(V)$ be a representation of $A_q$ such that there exists
a vector $\nu_\Lambda$ in $V$ satisfying
$$
\pi (b).\nu_\Lambda=\Lambda(b)\nu_\Lambda\,,~~~
\pi(\omega.b).\nu_\Lambda =V\,,~~~~\forall b\in {B}\,.
$$
Then $\pi$ is called a highest weight representation of $A_q$ with highest
weight $\Lambda$.

{\sf [Definition 2]}. Let $X$ be an m-dimensional manifold.
The current algebra
of $A$ is defined to be the algebra of continuous mappings
$$
A^X =Map(X,A)
$$
from $X$ to $A$.

We define in $A_q$ a topology as in \cite{jimbo}\cite{frenkel}.

{\sf [Definition 3]}. Let $X$ be an m-dimensional manifold.
The quantum current algebra
of $A_q$ is defined to be the quantum algebra of continuous mappings
$$
A_q^X =Map(X,A_q)
$$
from $X$ to $A_q$.

Let $X$ be a topological space and let $\mu$ be a finite positive measure.
We first investigate the continuous tensor product of quantum current
algebras and connect it to the highest weight representations. Set
$$
G_t =\left\{e^{ta}\vert a\in A_q\right\}\,,~~ t\in \Rb
$$
where $e^{ta}$ is defined by $e^{ta}=\displaystyle\sum_{n=0}^{\infty}
\frac{t^n}{n!}a^n$.
We call $G_t$ the $t$-exponential elements obtained from $A_q$ and the
fixed parameter $t\in \Rb$.

Set
$$
G_t^X =\left\{e^{t\sigma}\vert \sigma\in A_q^X\right\}\,, ~~t\in \Rb
$$
where $e^{t\sigma}$ is defined pointwise by $e^{t\sigma}(x)=e^{t\sigma(x)}$,
$\forall x\in X$. We call $G_t^X$ the t-exponential currents obtained from
the quantum current algebra of $A_q$. We give $G_t$ the topology induced
by the one defined in $A_q$.

Let $\{O_x\}_{x\in X}$ be an operator representation of $G_t$ acting on
some Hilbert space
${K}_X$,
$$
O_x:G_t\to u(K_X)
$$
with $u(K_X)$ the unitary operators in $K_X$.
Let $\{b_x\}_{x\in X}$ be an associated continuous $\mu$-measurable field of
one-cocycles,
$$
b_x:G_t\to K_X
$$
satisfying
\be\label{1}
b_x(g.g^{\prime})=b_x(g)+O_x(g)b_x(g^{\prime})\,,~~~~g,g^{\prime}\in G_t
\ee
such that for all $x$ in $X$, $b_x(G_t)$ is total in $K_X$. For an
example of such a cocycle we refer to section 3.
For all $f$ in $G_t^X$, we set
\be\label{2}
\tilde{O}(f)=\displaystyle\int_X^\oplus O_x[f(x)]d\mu(x)
\ee
\be\label{3}
\tilde{b}(f)=\displaystyle\int_X^\oplus b_x[f(x)]d\mu(x)
\ee
\be\label{4}
\tilde{K}=\displaystyle\int_X^\oplus K_X(x)\,d\mu(x)
\ee
with $K_X(x)$ the $x$-copy of $K_X$.
$\tilde{O}$ is a linear representation of $G_t$. The integrals have to be
understood as direct integrals, see e.g. \cite{gui}.

Let $\t{H}$ be the Fock space over $\t{K}$.
$$
\tilde{H}=Exp[\tilde{K}]\,.
$$
By construction the coherent states
$$
Exp[\tilde{k}]=
\sum_{n=0}^{\infty}\displaystyle\frac{\tilde{k}^{\otimes n}}{n!}\,,~~~
\tilde{k}\in \tilde{K}
$$
are total in $Exp[\tilde{K}]$, see e.g. \cite{gui}.
The representation $\tilde{u}=Exp[\tilde{O}]$
is then defined by the action on $Exp[\tilde{K}]$,
\be\label{5}
\tilde{u}(f).Exp[h]=e^{-\vert\vert\t{b}(f)\vert\vert^2/2
-<\t{O}(f).h,\t{b}(f)>}Exp[\t{O}(f).h+\t{b}(f)]\,,
\ee
where $h\in\tilde{K}$. $<~,~>$ denotes the scalar product in $\tilde{K}$
inherited from the one defined in the Hilbert space $K_x$, $x\in X$, i.e.
$$
<\tilde{k},\tilde{k^{\prime}}>
=\displaystyle\int_X^\oplus <k_x,k_x^{\prime}>d\mu(x)\,,~~~\forall
\tilde{k},~\tilde{k^{\prime}}\in \tilde{K}\,.
$$

For $\sigma\in A_q^x$, we define a representation $\pi$
of $\sigma$ in $Exp[\t{K}]$ by
\be\label{6}
\pi(\sigma)=\left\{\displaystyle\frac{d}{dt}\tilde{u}(f_t)\right\}_{t=0}\,,
\ee
where $f_t=e^{t\sigma}\in G^X_t$. Note that the derivative exists since
$\t{u}$ is an algebraic representation and $t\to f_t$ is also smooth
by construction.

{\sf [Proposition 1]}. $\pi$ acts as follows on the coherent vectors:
\be\label{7}
\ba{rcl}
\pi(\sigma).Exp[h]&=&
\left(-i\t{\phi}(\sigma)-<h,\t{\nu}_\sigma>\right)Exp[h]\\[4mm]
&&+{\displaystyle\sum_{n=0}^\infty}^\oplus
\displaystyle\frac{n+1}{\sqrt{(n+1)!}}[\t{\theta}(\sigma).h
+\t{\nu}_\sigma]\otimes h^{\otimes n}\,,~~~~\forall h\in\t{K}\,,
\ea
\ee
where
\begin{eqnarray}
\tilde{\theta}(\sigma)&=&
\left\{\displaystyle\frac{d}{dt}\t{O}(f_t)\right\}_{t=0}
{}~~~~\in End(\t{K})\\[4mm]
\tilde{\nu}_\sigma&=&\left\{\displaystyle\frac{d}{dt}\t{b}(f_t)\right\}_{t=0}
{}~~~~\in \t{K}\\[4mm]
\tilde{\phi}(\sigma)&=&
\left\{\displaystyle\frac{d}{dt}\t{\psi}(f_t)\right\}_{t=0}~~~~\in \Cb
\end{eqnarray}
$\pi$ is a unitary representation of $G^X_t$ on the Fock space
$\tilde{H}$. Moreover, if $\mu$ is non-atomic, $\pi$ is irreducible.

{\sf[Proof]}. This is a direct consequence of definitions of (\ref{5}) and
(\ref{6}) using methods in \cite{torresani}.  $\rule{2mm}{2mm}$

\section{The case of the q-deformed current algebra $sl_q(2,\Cb )^X$}

Now we take for $A_q$ the quantum algebra $A_q=sl_q(2,\Cb)$, as an example.
Let $e_0$, $f_0$ and $h_0$ be the generators of the algebra $A=sl(2,\Cb)$,
satisfying
\be\label{11}
[h_0,e_0]=e_0\,,~~~[h_0,f_0]=-f_0\,,~~~[e_0,f_0]=2h_0\,.
\ee
Let $e$, $f$ and $h$ be the corresponding
elements of the quantum algebra $sl_q(2,\Cb)$ (for the study of this algebra
see, e.g. \cite{jimbo}\cite{qg}) satisfying the commutation relations
\be\label{13}
[h,e]=e\,,~~~[h,f]=-f\,,~~~[e,f]=[2h]_q\,.
\ee
where $[x]_q=\displaystyle\frac{q^x-q^{-x}}{q-q^{-1}}$
and $q\in \Rb\backslash\{0\}$ is the deformation parameter.

The following mapping is a deformation map from
$sl(2,\Cb)$ to $sl_q(2,\Cb)$.
\be\label{12}
e=e_0\displaystyle\frac{[j-h]_q}{j-h}\,,~~~
f=f_0\displaystyle\frac{[j+h]_q}{j+h}\,,~~~h=h_0\,,
\ee
where $j$ is the center of $sl(2,\Cb)$.
It is straightforward to verify that $e,f,h$ satisfy the
$sl_q(2,\Cb)$ algebraic relations (\ref{13}) in terms of relations (\ref{11})
\cite{fg}\cite{zachos}.

Let $e_0(x),\,f_0(x),\,h_0(x)$ be the generators of the current algebra
$sl(2,\Cb)^X$ and
$e(x)$, $f(x)$, $h(x)$ the corresponding elements of the quantum current
algebra $sl_q(2,\Cb)^X$. Here and in the following $x$ runs over $X$.

{\sf [Proposition 2]}. The deformation map from
$sl(2,\Cb)^X$ to $sl_q(2,\Cb)^X$ is given by
\be\label{12x}
e(x)=e_0(x)\displaystyle\frac{[j(x)-h(x)]_q}{j(x)-h(x)}\,,~~~
f(x)=f_0(x)\displaystyle\frac{[j(x)+h(x)]_q}{j(x)+h(x)}\,,~~~h(x)=h_0(x)\,,
{}~~~~~\forall x\in X\,,
\ee
where $j(x)$ is the center of $sl(2,\Cb)^X$.

{\sf [Proof]}. The result is obvious from the map (\ref{12}) and the fact
that the current algebra is given the pointwise algebraic structure. Hence
the deformation mapping $A$ to $A_q$ directly gives rise to the
deformation mapping from $A^X$ to $A_q^X$. $\rule{2mm}{2mm}$

Therefore the map from a topological space $X$ to $A^X_q$ can be obtained from
the maps $X$ to $A^X$ and $A^X\to A_q^X$.

Now let us take $X=D$ with $D$ the unit disk of the complex plane
and consider the space $K=L^2_{hol}(D,\mu)$,
$L^2_{hol}$ being the space of square integrable holomorphic function on
$D$, $\mu$ being the Riemannian volume $\mu(dx)=r\,dr\,d\theta$,
$x=re^{i\theta}$, $r\in\Rb_+$, $\theta\in[0,2\pi]$.
The generators of $sl(2,\Cb)^D$ are
\be\label{14}
e_0(\xi)=\left(\begin{array}{cc}
0&\xi\\[3mm]
0&0\end{array}\right)\,,~~~
f_0(\xi)=\left(\begin{array}{cc}
0&0\\[3mm]
\xi&0\end{array}\right)\,,~~~
h_0(\xi)=\displaystyle\frac{1}{2}\left(\begin{array}{cc}
\xi&0\\[3mm]
0&-\xi\end{array}\right)\,,~~~\xi\in C_0^\infty (D)
\ee
with the center $j\to j(\xi)=\left(\begin{array}{cc}\xi&0\\[3mm]
0&\xi\end{array}\right)$, cfr. \cite{torresani}.

{}From the map (\ref{12}) we have the elements of $sl_q(2,\Cb)^D$,
\be\label{15}
e(\xi)=\displaystyle\frac{\sinh\gamma\xi}{\sinh\gamma}\left(\begin{array}{cc}
0&1\\[3mm]
0&0\end{array}\right)\,,~~~
f(\xi)=\displaystyle\frac{\sinh\gamma\xi}{\sinh\gamma}\left(\begin{array}{cc}
0&0\\[3mm]
1&0\end{array}\right)\,,~~~
h(\xi)=\displaystyle\frac{1}{2}\left(\begin{array}{cc}
\xi&0\\[3mm]
0&-\xi\end{array}\right)\,,
\ee
where $\gamma=\log q$ is the deformation parameter.

For $g=\left(\ba{cc}a&b\\[2mm]c&d\ea\right)\in G_t$,
$a,b,c,d\in C_0^\infty (D)$, the representation
$O_x:G\to u(K)$, $x\in D$, is defined by $O_x(g):f\in K\to O_x (g).f$,
\be\label{16}
[O_x(g).f](z)=[-cz+a]^{-2}f\left(\displaystyle\frac{dz+b}{-cz+a}\right)\,,
{}~~~\forall z\in D\,.
\ee
The $K$-valued one-cocycle is unique up to equivalence and is given by
\be\label{17}
[b_x(g)](z)=\displaystyle\frac{c}{a-cz}\,,~~~\forall z\in D\,.
\ee

For $f_t=e^{t\sigma}$, $\sigma=\left(\ba{cc}\alpha&\beta\\[2mm]
\gamma&\delta\ea\right)\in sl_q(2,\Cb)^D$, we have that the quantity (9) is
given by
\be\label{18}
\ba{rcl}
\t{\nu}_\sigma(z)&=&
\left\{\displaystyle\frac{d}{dt}\t{b}(f_t)(z)\right\}_{t=0}\\[4mm]
&=&\left\{\displaystyle\frac{d}{dt}\displaystyle\int_D^\oplus
b_x[f_t(x)](z)d\mu(x)\right\}_{t=0}\\[4mm]
&=&\displaystyle\int_D^\oplus \gamma(x)d\mu(x)\lambda_0(z)\,,~~~z\in D
\ea
\ee
with $\lambda_0$ the unit function on $D$ defined by $\lambda_0(z)=1$.
Correspondingly $\tilde{\phi}$, occurring in (10), is here given by
\be\label{19}
\tilde{\phi}(\sigma)=
\left\{\displaystyle\frac{d}{dt}\t{\psi}(f_t)\right\}_{t=0}
=\displaystyle\int_D\alpha(x)d\mu(x)\,.
\ee

{\sf [Proposition 3].} Let $\pi$ defined by (\ref{7}) as a representation of
$sl_q(2,\Cb)^D$ in $Exp[\tilde{K}]$,
$\tilde{K}\equiv\int_D^\oplus K_D(x)\,d\mu(x)$.
Then $\pi$ is a highest weight representation of $sl_q(2,\Cb)^D$, and its
highest weight is completely determined by the Riemannian volume
$\mu$ on $D$.

{\sf [Proof].} We have only to check that $\pi$ admits a highest weight vector.
Let $\Omega=Exp[0]$, where $0$ is the zero element in $\t{K}$.
{}From (\ref{7}), (\ref{15}), (\ref{18}) and (\ref{19})
we have
\be\label{20}
\ba{rcl}
\pi(e(\xi))&=&0\\[3mm]
\pi(f(\xi))&=&
\displaystyle\int_D^\oplus\displaystyle\frac{\sinh\gamma\xi(x)}{\sinh\gamma}
d\mu(x)\lambda_0 \\[3mm]
\pi(h(\xi))&=&
\left[-i\displaystyle\int_D\displaystyle\frac{\sinh\gamma\xi(x)}
{\sinh\gamma}d\mu(x)\right].\Omega
\ea
\ee
Therefore from the definition 1, $\Omega$ is a highest weight vector for $\pi$.
The corresponding highest weight of the quantum current algebra
$sl_q(2,\Cb)^D$ is completely determined by the $\mu$ measure
on $D$. $\rule{2mm}{2mm}$

We remark that when the deformation parameter $\gamma$ approaches zero,
formulae (\ref{20}) become the ones of $sl(2,\Cb)^D$ (in the sense that
the formal replacement of $\gamma$ by zero in (\ref{20}) yields the
formulae in \cite{torresani}). Hence our
representation is a q-extension of the representation of $sl(2,\Cb)^D$
obtained in \cite{torresani}.

Summarizing, we have discussed the representations of quantum current
algebras defined by mappings from manifolds into quantum
algebras. This is a continuous tensor product representation.
It is shown that for the quantum current algebra $sl_q(2,\Cb)^D$, the
continuous tensor product representation is also a highest weight
representation.

{\sf Remark:}
We have taken the Hilbert space to be real in this paper. Our construction
can be extended to the case where the Hilbert space
is complex, by changing the definition (\ref{6}) by a factor corresponding
to a nonvanishing 2-cocycle \cite{abook}\cite{torresani}.

\vspace{2.5ex}
ACKNOWLEDGEMENTS: We are very grateful to Professor Torr\'esani for very
helpful discussions in connection with his work which has been the
basis for our own work. We thank A.v. Humboldt foundation for the financial
support given to the second named author.

\end{document}